\documentclass[12pt]{amsart}

\usepackage{amsmath,amsfonts,amssymb,amscd,verbatim,delarray,verbatim}
\newcommand{\F}{{\mathbb F}}

\newcommand{\Q}{{\mathbb Q}}

\begin{document}
\numberwithin{equation}{section}

\newtheorem{theorem}{Theorem}[section]
\newtheorem{lemma}[theorem]{Lemma}
\newtheorem{prop}[theorem]{Proposition}
\newtheorem{proposition}[theorem]{Proposition}
\newtheorem{corollary}[theorem]{Corollary}
\newtheorem{corol}[theorem]{Corollary}
\newtheorem{conj}[theorem]{Conjecture}

\theoremstyle{definition}
\newtheorem{defn}[theorem]{Definition}
\newtheorem{example}[theorem]{Example}
\newtheorem{examples}[theorem]{Examples}
\newtheorem{remarks}[theorem]{Remarks}
\newtheorem{remark}[theorem]{Remark}
\newtheorem{algorithm}[theorem]{Algorithm}
\newtheorem{question}[theorem]{Question}
\newtheorem{problem}[theorem]{Problem}
\newtheorem{subsec}[theorem]{}


\def\toeq{{\stackrel{\sim}{\longrightarrow}}}
\def\into{{\hookrightarrow}}


\def\alp{{\alpha}}  \def\bet{{\beta}} \def\gam{{\gamma}}
 \def\del{{\delta}}
\def\eps{{\varepsilon}}
\def\kap{{\kappa}}                   \def\Chi{\text{X}}
\def\lam{{\lambda}}
 \def\sig{{\sigma}}  \def\vphi{{\varphi}} \def\om{{\omega}}
\def\Gam{{\Gamma}}   \def\Del{{\Delta}}
\def\Sig{{\Sigma}}   \def\Om{{\Omega}}
\def\ups{{\upsilon}}


\def\F{{\mathbb{F}}}
\def\BF{{\mathbb{F}}}
\def\BN{{\mathbb{N}}}
\def\Q{{\mathbb{Q}}}
\def\Ql{{\overline{\Q }_{\ell }}}
\def\CC{{\mathbb{C}}}
\def\R{{\mathbb R}}
\def\V{{\mathbf V}}
\def\D{{\mathbf D}}
\def\BZ{{\mathbb Z}}

\def\XX{\mathbf{X}^*}
\def\xx{\mathbf{X}_*}

\def\AA{\Bbb A}
\def\BA{\mathbb A}
\def\HH{\mathbb H}
\def\PP{\Bbb P}

\def\Gm{{{\mathbb G}_{\textrm{m}}}}
\def\Gmk{{{\mathbb G}_{\textrm m,k}}}
\def\GmL{{\mathbb G_{{\textrm m},L}}}
\def\Ga{{{\mathbb G}_a}}

\def\Fb{{\overline{\F }}}
\def\Kb{{\overline K}}
\def\Yb{{\overline Y}}
\def\Xb{{\overline X}}
\def\Tb{{\overline T}}
\def\Bb{{\overline B}}
\def\Gb{{\bar{G}}}
\def\Ub{{\overline U}}
\def\Vb{{\overline V}}
\def\Hb{{\bar{H}}}
\def\kb{{\bar{k}}}

\def\Th{{\hat T}}
\def\Bh{{\hat B}}
\def\Gh{{\hat G}}

\def\cF{{\mathfrak{F}}}
\def\cC{{\mathcal C}}
\def\cU{{\mathcal U}}

\def\Xt{{\widetilde X}}
\def\Gt{{\widetilde G}}

\def\gg{{\mathfrak g}}
\def\hh{{\mathfrak h}}
\def\lie{\mathfrak a}

\def\min{^{-1}}

\def\textrm#1{\text{\textnormal{#1}}}

\def\GL{\textrm{GL}}            \def\Stab{\textrm{Stab}}
\def\Gal{\textrm{Gal}}          \def\Aut{\textrm{Aut\,}}
\def\Lie{\textrm{Lie\,}}        \def\Ext{\textrm{Ext}}
\def\PSL{\textrm{PSL}}          \def\SL{\textrm{SL}}
\def\loc{\textrm{loc}}
\def\coker{\textrm{coker\,}}    \def\Hom{\textrm{Hom}}
\def\im{\textrm{im\,}}           \def\int{\textrm{int}}
\def\inv{\textrm{inv}}           \def\can{\textrm{can}}
\def\id{\textrm{id}}              \def\char{\textrm{char}}
\def\Cl{\textrm{Cl}}
\def\Sz{\textrm{Sz}}
\def\ad{\textrm{ad\,}}

\def\tors{_{\textrm{tors}}}      \def\tor{^{\textrm{tor}}}
\def\red{^{\textrm{red}}}         \def\nt{^{\textrm{ssu}}}

\def\sss{^{\textrm{ss}}}          \def\uu{^{\textrm{u}}}
\def\mm{^{\textrm{m}}}
\def\tm{^\times}                  \def\mult{^{\textrm{mult}}}

\def\uss{^{\textrm{ssu}}}         \def\ssu{^{\textrm{ssu}}}
\def\comp{_{\textrm{c}}}
\def\ab{_{\textrm{ab}}}

\def\et{_{\textrm{\'et}}}
\def\nr{_{\textrm{nr}}}

\def\nil{_{\textrm{nil}}}
\def\sol{_{\textrm{sol}}}
\def\End{\textrm{End\,}}

\def\til{\;\widetilde{}\;}


\title[characterization of radical]
{{\bf Engel-like   characterization of radicals\\
in finite dimensional Lie algebras\\ and  finite  groups }}
\author[Bandman, Borovoi, Grunewald, Kunyavski\u\i , Plotkin]{
Tatiana Bandman, Mikhail Borovoi, Fritz Grunewald,\\
Boris Kunyavski\u\i ,  and Eugene Plotkin }
\address{Bandman, Kunyavski\u\i \ and Plotkin: Department of
Mathematics, Bar-Ilan University, 52900 Ramat Gan,
ISRAEL}
\email{bandman@macs.biu.ac.il , kunyav@macs.biu.ac.il,
plotkin@macs.biu.ac.il}
\address{Borovoi: Raymond and Beverly Sackler School of Mathematical Sciences, Tel Aviv University,
69978 Tel Aviv, ISRAEL} \email{borovoi@post.tau.ac.il }
\address{Grunewald: Mathematisches Institut der Universit\"at
Heinrich Heine D\"usseldorf, Universit\"atsstr. 1, 40225
D\"usseldorf, GERMANY}
\email{grunewald@math.uni-duesseldorf.de}

\begin{abstract}
A classical theorem of R. Baer describes the nilpotent radical of
a finite group $G$ as the set of all Engel elements, i.e. elements
$y\in G$ such that for any $x\in G$ the $n$th commutator
$[x,y,\dots ,y]$ equals 1 for $n$ big enough. We obtain a
characterization of the solvable radical of a finite dimensional
Lie algebra defined over a field of characteristic zero in similar
terms. We suggest a conjectural description of the solvable
radical of a finite group as the set of Engel-like elements and
reduce this conjecture to the case of a finite simple group.
\end{abstract}

\maketitle

\thispagestyle{empty} \vspace{1.0cm} 

\tableofcontents

\section{Introduction} \label{sec:intro}

In the present paper we study certain Engel-like characterizations of the
solvable radical of finite dimensional Lie algebras and finite
groups. Such a characterization for the nilpotent radical of a
finite group is given by R.~Baer \cite{Ba}, see also \cite{H}, \cite{Ro}:

\begin{theorem}\label{th:br}
The nilpotent radical of a finite group $G$ coincides with the
collection of all Engel elements of $G$.
\end{theorem}

This theorem of Baer generalizes Zorn's theorem \cite{Zo} which
gives a characterization of finite nilpotent groups in terms of
special two-variable Engel identities.

Our goal is to establish an analog of Baer's theorem for the solvable
radical of a finite group. As a first step, we consider a similar problem
for finite dimensional Lie algebras.
Towards this end, we introduce some
sequences of words which provide the solvability property. These
sequences play the role of the Engel sequences that work in the
nilpotent case.

All groups in the paper are assumed to be finite, all Lie algebras
are finite dimensional. Observe, however, that Baer's theorem is
valid for many classes of infinite groups satisfying some
finiteness conditions. For example, Baer himself proved it for
noetherian groups, a similar result is valid for radical groups
\cite{Plo1}, \cite{Plo2}, linear groups \cite{Pla}, and PI-groups
\cite{Plo3} (see also \cite{Pr}, \cite{To}).

The paper is organized as follows. In Section 2 all main
definitions and statements of the problems are given. Section 3
deals with the case of Lie algebras. In Section 4 we prove a
reduction theorem for the case of groups.



\noindent {\it Acknowledgements}. Bandman, Kunyavski\u\i , and
Plotkin were partially supported by the Ministry of Absorption
(Israel), the Israeli Science Foundation founded by the Israeli
Academy of Sciences --- Center of Excellence Program, and the
Minerva Foundation through the Emmy Noether Research Institute of
Mathematics. Borovoi, Kunyavski\u\i \ and Plotkin were also
supported by the RTN network HPRN-CT-2002-00287, and Kunyavski\u\i
\ and Plotkin by INTAS 00-566.

We are very grateful to A.~Merkurjev, B.~Plotkin, A.~Premet,
M.~Sapir, and E.~Vinberg for useful discussions and
correspondence.

\section{Definitions and main problems} \label{sec:def}

We follow the terminology introduced in \cite{Plo3}.

Let $F_2=F(x,y)$ be the free two generator group, $W_2=W(x,y)$
the free two generator Lie algebra.

\begin{defn}
We say that a sequence $\overrightarrow{u} = u_1, u_2, u_3, \dots,
u_n, \dots$ of elements from $F_2$ is {\bf correct} if the
following conditions hold:

(i) $u_n(a,1)=1$ and $u_n(1,g)=1$ for all sufficiently big $n$,
every group $G$, and all elements $a,g\in G$;

(ii) if $a,g$ are elements of $G$ such that $u_n(a,g)=1$, then for
every $m>n$ we have $u_m(a,g)=1$.
\end{defn}

Thus, if the identity $u_n(x,y)\equiv 1$ is satisfied in $G$ then
for every $m>n$ the identity  $u_m(x,y)\equiv 1$ also holds in
$G$.

A similar definition can be given for Lie algebras.

\begin{defn}
For every correct sequence $\overrightarrow{u}$ define the class
of groups (resp. Lie algebras) $\Theta=\Theta
(\overrightarrow{u})$ by the rule: a group (resp. Lie algebra) $G$
belongs to $\Theta$ if and only if there is $n$ such that the identity
$u_n(x,y)\equiv 1$ (resp. $u_n(x,y)\equiv 0$) holds in $G$.
\end{defn}

\begin{defn}
For every group (resp. Lie algebra) $G$ denote by
$G(\overrightarrow{u})$ the subset of $G$ defined by the rule:
$g\in G(\overrightarrow{u})$ if and only if for every $a\in G$ there exists
$n=n(a,g)$ such that $u_n(a,g)=1$ (resp. $0$). Elements of
$G(\overrightarrow{u})$ are viewed as Engel elements with respect
to the given correct sequence $\overrightarrow{u}$. We call these
elements {\bf $\overrightarrow{u}$-Engel-like} or, for brevity,
{\bf $\overrightarrow{u}$-Engel elements}.
\end{defn}

\begin{examples}

${}$

\begin{enumerate}

\item If $\overrightarrow{u}=\overrightarrow{e}=e_1, e_2,
e_3,\ldots, $ where
$$
e_1(x,y)=[x,y]=xyx^{-1}y^{-1},\ldots, e_n(x,y)=[e_{n-1}(x,y),y],\ldots,
$$
then $\Theta(\overrightarrow{e})$ is the class of all Engel
groups. In the case of finite groups the class
$\Theta(\overrightarrow{e})$ coincides with the class of
finite nilpotent groups \cite{Zo}.

Clearly, $\overrightarrow{e}$-Engel elements of any group $G$ are
none other than usual Engel elements. If $G$ is finite, the set
$G(\overrightarrow{e})$ coincides with the nilpotent radical of
$G$ \cite{Ba}.

\item If $\overrightarrow{u}$ is determined by the following correct
sequence of words:
$$
u_1=xy^{-1}, u_2=u_1(xy,yx)=[x,y], \ldots, u_n=u_{n-1}(xy,yx),\dots ,
$$
then for finite groups the class $\Theta(\overrightarrow{u})$
coincides with the class of all finite nilpotent-by-two groups
\cite{BP}.

\item The correct sequence $\overrightarrow{u}$, where
$$
u_1(x,y)=x^{-2}y^{-1}x,  \ldots,  u_{n}(x,y)=[\, x\, u_{n-1}
(x,y)\, x^{-1},\, y\, u_{n-1}(x,y)\, y^{-1}\, ],\dots ,
$$
determines the class $\Theta(\overrightarrow{u})$ which coincides with
the class of finite solvable groups \cite{BGGKPP}.

\item  The correct sequence $\overrightarrow{s}$, where
$$
s_1(x,y)=x, \ldots, s_n(x,y)=[s_{n-1}(x,y)^{-y},
s_{n-1}(x,y)],\ldots
$$
determines the class $\Theta(\overrightarrow{s})$ which also coincides
with the class of finite solvable groups \cite{BWW}.

\item The correct sequence $\overrightarrow{w}$, where
$$
w_1(x,y)=[x,y],\ldots,w_n(x,y)=[[w_{n-1},x],[w_{n-1},y]], \ldots
$$
and $[\ ,\ ]$ stands for the Lie bracket in a Lie algebra, determines the class
$\Theta(\overrightarrow{w})$ of finite dimensional solvable Lie
algebras over an infinite field $k$, $\char k\neq 2,3,5$
\cite{GKNP}.
\end{enumerate}
\end{examples}

It is easy to see that if $g$ is a $\overrightarrow{u}$-Engel
element in $G$ then $g$ is a $\overrightarrow{u}$-Engel element
in every subgroup $H$ containing $g$. If $H$ is a normal
subgroup in $G$ and $g$ is a $\overrightarrow{u}$-Engel element in $G$,
then $\bar g=gH$ is a
$\overrightarrow{u}$-Engel element in $\bar G=G/H$.

The following natural problem arises:

\begin{problem}\label{ref:solvgroups}
Describe the class of sequences $\overrightarrow{u}$ such that
$\Theta(\overrightarrow{u})$ is the class of finite solvable
groups.
\end{problem}

Denote by $F_2^{(n)}(x,y)$ the $n$-th term of the derived series
of the free group $F_2(x,y)$.

\begin{proposition}
Let $\overrightarrow{u}$ be a correct sequence. Then
$\Theta(\overrightarrow{u})$ coincides with the class of finite
solvable groups if and only if

(i) for every $n$ there exists $k=k(n)$ such that $u_k$ belongs to
$F_2^{(n)}(x,y)$;

(ii) there is no $n$ such that $u_n(x,y)=1$ is an identity in one of the
following groups: {\rm (1)} $G=\PSL (2,\BF_q)$ where $q\ge 4$ is a
prime power, {\rm (2)} $G=\Sz (2^m)$, $m\in\BN,\, m\ge 3$
and odd, {\rm (3)} G=$\PSL (3,\BF_3)$.

Here $\PSL (m,\mathbb F_q)$ denotes the projective special linear
group of degree $m$ over $\mathbb F_q$. For $q=2^m$ we denote by
$\Sz (q)$ the Suzuki group (the twisted form of ${}^2\!B_2$, see
\cite[XI.3]{HB}).
\end{proposition}

\begin{proof} The proof repeats the proof of Theorem 2.1 from
\cite{BGGKPP} and is based on the list of minimal simple
non-solvable groups \cite{Th}.
\end{proof}

\begin{remark} Even for the case of the group $G=\PSL (2,\BF_q)$
the basis of identities is known only for small fields $\BF_q$ and
not known in general \cite{CMS}, \cite{So}. The bases look highly
complicated and do not provide any hint to check explicitly if a
particular identity follows from the basis.
\end{remark}

\begin{problem} \label{ref:solalgebras}
Describe the class of sequences $\overrightarrow{u}$ such that
$\Theta(\overrightarrow{u})$ is the class of finite dimensional
solvable Lie algebras.
\end{problem}

Denote by $W_2^{(n)}(x,y)$ the $n$-th term of the derived series
of the free Lie algebra $W_2(x,y)$.

\begin{proposition} Let $\overrightarrow{u}$ be a correct sequence.
Then $\Theta(\overrightarrow{u})$ coincides with the class of
finite dimensional Lie algebras over an infinite field $k$,
$\char\ k\neq 2,3,5$ if and only if

(i) for every $n$ there exists $k=k(n)$ such that $u_k$
belongs to $W_2^{(n)}(x,y)$;

(ii) there is no $n$ such that $u_n(x,y)=0$ is an identity in the simple
Lie algebra $\text{\rm{sl}}_2$.
\end{proposition}

\begin{proof} The proof immediately follows from the fact that if
the characteristic of $k$ differs from 2,3,5, then every simple
Lie $k$-algebra contains a subalgebra which has the algebra $\text{\rm{sl}}_2$
as a quotient (see \cite{GKNP} and references therein).
\end{proof}

\begin{remark}
The assumption $\char k \neq 3,5$ is technical and can probably be
dropped.
\end{remark}

\begin{remark} The situation with bases of identities for the minimal
simple non-solvable Lie algebras is different from the group case.
A basis of identities of the algebra $\text{\rm{sl}}_2$ over a
field of characteristic zero is known and consists of two
identities \cite{Ra}, \cite{Bah} (see \cite{MK} for the case where
the ground field is finite, and \cite{Ko} and references therein
for the case where the ground field is infinite of characteristic
$p$). However, it is difficult to verify whether a particular
$u_n$ is an identity in $\text{\rm{sl}}_2$ using this basis.
\end{remark}


\begin{conj} \label{conj:radgr}
There is a sequence $\overrightarrow{u}=\overrightarrow{u}(x,y)$
such that for every finite group $G$ the solvable radical of $G$
coincides with $G(\overrightarrow{u})$.
\end{conj}

\begin{conj}\label{conj:radlie}
There is a sequence $\overrightarrow{u}=\overrightarrow{u}(x,y)$
such that for every finite dimensional Lie algebra $L$ the
solvable radical of $L$ coincides with $L(\overrightarrow{u})$.
\end{conj}

Of course, it is highly desirable not only to prove the existence
of such a $\overrightarrow{u}$, but to exhibit an explicit
sequence in each of the above cases (groups and Lie algebras).
We give a partial positive answer to Conjecture \ref{conj:radlie}
in Section \ref{sec:lie} (see Theorems \ref{th:rad} and
\ref{th:rad-w}), and reduce Conjecture \ref{conj:radgr} to the
case of simple groups in Section \ref{sec:gr} (see Theorem
\ref{th:red}).


\section{Case of Lie algebras}\label{sec:lie}

Let $L$ be a finite dimensional Lie algebra over a field $k$.
Denote by $[ , ]$ the Lie operation. For $t\in L$ the linear
operator $\ad t\colon L\to L$ is defined by $(\ad t)x = [t,x]$.

By the {\em solvable radical} of $L$ we mean the largest solvable ideal $R$ of $L$
(Bourbaki \cite{Bou} and Jacobson \cite{J} call $R$ the radical of $L$).
By the {\em nilpotent radical} of $L$ we mean the largest nilpotent ideal $N$ of $L$
(Jacobson \cite{J} calls $N$ nil radical, and Bourbaki  \cite{Bou} calls it just
the largest nilpotent ideal).

Define the sequence $\overrightarrow{e}$ by
$e_1(x,y)=[x,y]$ and, by induction,  $e_{n+1}(x,y)=[e_n(x,y),y]$.
Then $e_{n}=(-\ad y)^n x$.

\begin{defn} \label{def:engel}
An element $y\in L$ is called an {\em Engel} element if it is
$\overrightarrow{e}$-Engel, i.e. for every $x\in L$ there exists $n$
such that $e_n(x,y)=0$ (i.e. $(\ad y)^n x=0$).
\end{defn}

The following proposition is well known.

\begin{proposition}\label{prop:engel}
An element $y\in L$ is Engel if and only if $\ad y$ is nilpotent
(i.e. there exists $n$ such that $(\ad y)^n=0$).
\end{proposition}

We need a lemma.

\begin{lemma} \label{lem:dim}
Let $V$ be a vector space of dimension $d$ over a field $k$, let
$x\in V$, and let $A\colon V\to V$ be a linear map. If $A^mx=0$
for some $m$, then $A^dx=0$.
\end{lemma}

\begin{proof}[Proof of the lemma]
By Fitting's lemma (see \cite[Sect. II.4]{J}), $V=V_0\oplus V_1$,
where $V_0$ and $V_1$ are invariant subspaces of $A$,
the restriction $A_0$ of $A$ to $V_0$ is nilpotent,
and the restriction $A_1$ of $A$ to $V_1$ is invertible.
Write $x=(x_0,x_1)$, where $x_i\in V_i$ $(i=0,1)$. We have
$A_1^m x_1=0$, hence $x_1=0$ (because $A_1$ is invertible), i.e.
$x\in V_0$. Since $A_0$ is nilpotent, we have $A_0^{d_0}=0$ where
$d_0=\dim V_0$. We thus have $A^{d_0}x=0$ (because $x\in V_0$).
Since $d\geq d_0$, we conclude that $A^d x=0$.
\end{proof}

\begin{proof}[Proof of the proposition]
If $\ad y$ is nilpotent, then clearly $y$ is Engel.
Conversely, if $y$ is Engel, then by Lemma \ref{lem:dim}
$(\ad y)^d x=0$ for any $x\in L$, where $d=\dim L$.
Hence $(\ad y)^d=0$ and $\ad y$ is nilpotent.
\end{proof}

Define the sequence $\overrightarrow{v}$ by $v_1(x,y)=x$ and, by
induction, $v_{n+1}(x,y)=[v_n(x,y),[x,y]]$.
Then $v_{n+1}(x,y)=(-\ad [x,y])^n x =e_n(x,[x,y])$.

\begin{theorem}\label{th:ch}
Let $L$ be a finite dimensional Lie algebra over a field $k$ of
characteristic zero. Then $L$ is solvable if and only if for some
$n$ the identity $v_n(x,y)\equiv 0$ holds in $L$.
\end{theorem}

\begin{proof} If $L$ is solvable, then $L'=[L,L]$ is nilpotent
\cite[Cor. II.7.1]{J}. Hence every pair $z,t$ of elements of
$L'$ satisfies the identity  $(\ad t)^m z=0$, where $m=\dim L'$.
On putting $z=[x,[x,y]]$, $t=[x,y]$, we get $v_{m+2}(x,y)=0$.

In the opposite direction, we mimic the proof of
\cite[Thm. 3.1]{GKNP}. Namely, we first reduce to the case when $k$
is algebraically closed. This reduction immediately follows from
the following easy lemma (cf. \cite{GKNP}).

\begin{lemma} \label{lem:cl}
Let $k$ be an infinite field and let $P(x_1,\dots, x_n)$ be a polynomial in $n$
variables over $k$.
If $P(x_1,\dots, x_n)=0$ for all $x_1,\dots, x_n\in k$, then
for any field extension $K/k$ we have
$P(x_1,\dots, x_n)=0$ for all $x_1,\dots, x_n\in K$.
\end{lemma}
\begin{proof}
Since $k$ is infinite, by \cite[Cor.~IV.1.7]{L}
$P$ is the zero polynomial, and the lemma follows.
\end{proof}


\begin{remark} A similar statement is valid for the elements of
free associative and free Lie algebras, i.e. for non-commutative
associative and Lie polynomials.
\end{remark}


Thus we now assume $k$ algebraically closed and suppose that $L$
satisfies $v_n(x,y)\equiv 0$ and is not solvable. Denote by $R$
the solvable radical of $L$. The algebra $L/R$ is not zero and
semisimple and also satisfies the identity $v_n(x,y)\equiv 0$. It
contains a subalgebra $S$ isomorphic to $\text{sl}_2$ which also must
satisfy the same identity. Being isomorphic to $\text{sl}_2$,
the Lie algebra $S$ has a basis
$\{e_+, e_-, h\}$ such that $[h,e_+]=2e_+,\ [h,e_-]=-2e_-,\ [e_+,e_-]=h$.
 Take $x=e_+$, $y=e_-$.
 We have $[x,y]=h$, hence
$v_2(x,y)=[e_+,h]=-2e_+$, and for any $n$
we get $v_n(x,y)=(-2)^{n-1}e_+\neq 0$, contradiction. The
theorem is proved.
\end{proof}

It turns out that the same sequence $v_n$ allows one to describe
the solvable radical of $L$.


\begin{theorem} \label{th:rad}
Let $L$ be a finite dimensional Lie algebra over a field $k$ of
characteristic zero. Then its solvable radical $R$ coincides with
the set of all $\overrightarrow{v}$-Engel elements of $L$.
\end{theorem}

\begin{proof}
Let us first prove that every element $y$ of $R$ is
$\overrightarrow{v}$-Engel. Since $\text{char}(k)=0$, by
\cite[Thm. II.7.13]{J}, $[L,R]$ is a nilpotent ideal. Hence every
element $t$ of $[L,R]$ is Engel in $L$ \cite[Thm. II.3.3]{J}, i.e. for any
$z\in L$ there exists $n$ such that $(\ad t)^n z=0$.
Let $x$ be an arbitrary element of $L$. The
element $[x,y]$ is Engel in $L$, hence on setting $t=[x,y]$, $z=x$, we
get $v_n(x,y)=0$ for some $n$.

Conversely, let us prove that every
$\overrightarrow{v}$-Engel element lies in the solvable radical.
Let us first consider the case where $k$ is algebraically closed.
It is enough to prove that there are no nonzero $\overrightarrow{v}$-Engel
elements in a simple $k$-algebra $L$.

Let $L$ be a simple Lie algebra over $k$. Choose a Cartan
subalgebra $H\subset L$ and a Borel subalgebra $B\supset H$. Let
$$
L=H\oplus\bigoplus_{\alpha\in\Phi}L_\alpha
$$
be the root decomposition. Here $\Phi=\Phi (L,H)$ is a reduced
irreducible root system, and all the spaces $L_\alpha$ are
1-dimensional. Let $\Phi^+$ be the set of positive roots
determined by $B$. We have
$$
L=L^-\oplus H\oplus L^+, \text{ where }
L^+=\bigoplus_{\alpha\in\Phi^+}L_\alpha,\
L^-=\bigoplus_{\alpha\in\Phi^+}L_{-\alpha}\;.
$$
We write an element $y\in L$ in the form $y=y^- +h+ y^+$, where
$y^-\in L^-,\ h\in H,\ y^+\in L^+$. For every $\alpha\in \Phi$
choose  a nonzero element $e_\alpha\in L_\alpha$. The set
$\{e_{-\alpha}:\ \alpha\in\Phi^+\}$ is a basis in $L^-$. Write
$y^-=\sum_{\alpha\in\Phi^+}c_{-\alpha}(y)e_{-\alpha}$, where
$c_{-\alpha}(y)\in k$.

Let $y\in L$, $y\neq 0$, $y=y^- +h+ y^+$.
Since $\Phi$ is irreducible, it has a unique maximal root;
we denote it by $\gamma$.
We may assume that $H$ and $B$ are chosen so that the
coefficient  $c_{-\gamma}(y)$ of $y^-$ at $e_{-\gamma}$ is nonzero.

Indeed, let $\Aut L$ denote the algebraic group of automorphisms
of $L$, and let $G$ be the connected component of $\Aut L$. It
suffices to prove that there exists $g\in G(k)$ such that
$c_{-\gamma}(gy)\neq 0$. Assume the contrary, i.e. that
$c_{-\gamma}(gy)= 0$ for any $g$. Set $V=\{z\in
L\;|\;c_{-\gamma}(z)=0\}$, it is a subspace of $L$ of codimension 1.
Let $W(y)$ denote the vector space generated by the elements of
the form $gy$ for $g\in G(k)$, it is a nonzero $G$-invariant
subspace of $L$. Since by assumption $W(y)\subset V$, we conclude
that $W(y)\neq L$. This leads to a contradiction, because the
representation of $G$ in $L$ is irreducible for a simple Lie
algebra $L$.

So let $y\in L,\ y\neq 0$, and $c:=c_{-\gamma}(y)\neq 0$.
We shall find $x$ such that
$v_n(x,y)\neq 0$ for any $n$. Take $x=e_\gamma$. Denote
$h_{\gam}=[e_{\gam},e_{-\gam}]$, it is a nonzero element of $H$.
We have $[h_\gamma,e_\gamma]=a e_\gamma$ with $a\neq 0$.
Since $\gamma$ is the maximal root, we have
$[x,y]=ch_{\gam}+y_1^+$ where $y_1^+\in L^+$.
Hence we have
$[x,[x,y]]=[e_{\gamma},ch_{\gam}+y_1^+]=-cae_{\gam}$, and, by induction,
$v_n(x,y)=(-ca)^{n-1} e_{\gam}$.
Since $a\neq 0$ and $c\neq 0$, we see that
$v_n(x,y)\neq 0$ for all $n$. Thus there are no nonzero
$\overrightarrow{v}$-Engel elements in $L$. This proves our
theorem in the case when $k$ is algebraically closed.

\bigskip

To reduce the general case to the case of the algebraically closed
ground field, we shall show that if $y$ is a
$\overrightarrow{v}$-Engel element of a finite dimensional Lie
algebra $L$, then one can choose $n$ in the condition $v_n(x,y)=0$
to be independent of $x$.
Actually we shall show that one can choose $n$ to be independent of $x$ and $y$.

\begin{lemma} \label{cor:v}
Assume that $v_n(x,y)=0$ for some $x,y\in L$ and for some natural
$n$. Then $v_{d+1}(x,y)=0$, where $d=\dim L$.
\end{lemma}
\begin{proof}
We apply Lemma \ref{lem:dim} to the linear operator
$A=-\ad [x,y]\colon z\mapsto [z,[x,y]]$ in the linear space $L$.
\end{proof}

\def\Rb{{\overline{R}}}
\def\Lb{{\overline{L}}}

We can now complete the proof of Theorem \ref{th:rad}. Let $y\in
L$ be a $\overrightarrow{v}$-Engel element of $L$. This means that
for any $x\in L$ there exists $n$ such that $v_n(x,y)=0$. By
Lemma \ref{cor:v}, we then get $v_{d+1}(x,y)=0$ for any $x\in
L$. This is a polynomial identity in $x$. Since the field $k$ is
infinite, it follows from Lemma \ref{lem:cl}
that $v_{d+1}(\bar x,y)=0$ for any $\bar x\in\Lb$,
where $\Lb=L\otimes\bar{k}$. 
In other words, $y$ is a $\overrightarrow{v}$-Engel element of $\Lb$.
Since our theorem is already proved over an algebraically closed
field, we see that  $y\in \Rb$, where $\Rb$ is the solvable radical of $\Lb$.
But $\Rb=R\otimes_k\bar k$ \cite[Ch.~I, \S5, n${}^\circ$6]{Bou},  hence $y\in R$.
\end{proof}

\begin{remark}
The sequence $\overrightarrow{v}$ is adjusted to the case of Lie
algebras over a field of characteristic zero. Indeed, the key
point in the proof of Theorems \ref{th:ch} and \ref{th:rad} was
the fact that if $L$ is solvable then $[L,L]$ is nilpotent. This
is no longer true in  positive characteristic. Here is an explicit
counter-example to the corresponding statements
in positive characteristic, based on the fact that $[L,R]$ is
no longer contained in the nilpotent radical of $L$.
\end{remark}


\begin{example} \label{ex:J}
We use an example given in \cite[II.7, pp.~52--53]{J}.
Let $k$ be a field of characteristic $p>0$, and $L$ be a
vector $k$-space of dimension $p+2$. Denote by $\{e,f,e_1,\dots ,e_p\}$
a basis of $L$ and define a structure of Lie algebra by the following
multiplication table:

\begin{equation} \label{eq:J}
\begin{array}{lll}
[e,f]=e, & [e,e_i]=e_{i+1} \quad (1\leq i\leq p-1), &  [e,e_p]=e_1, \\
& [f,e_i]=(i-1)e_i \quad (1\leq i\leq p),& [e_i,e_j]=0 \quad (1\leq i,j\leq p).
\end{array}
\end{equation}

From formulas (\ref{eq:J}) it follows immediately that the subalgebra
$M=\left<e_1,\dots ,e_p\right>$ is an abelian ideal of $L$ and the quotient
$S=L/M$ is a two-dimensional solvable algebra, hence $L$ is solvable (see
\cite[{\it loc.~cit.}]{J}). We shall show that $L$ does not satisfy any of
identities $v_n(x,y)\equiv 0$. Indeed, take $x=f+e_1$, $y=e+e_2$. We have
$$
t:=[x,y]=-e,
v_1(x,y)=x,
v_2(x,y)=[x,t]=e+e_2,
v_3(x,y)=e_3,\dots ,
$$
and, by induction, $v_p=e_p$, $v_{p+1}=e_1$, \dots ,
$v_{mp+r}=e_r$, \dots . Thus for all $n$ we have $v_n(x,y)\neq 0$,
and Theorem \ref{th:ch}  fails in characteristic $p$.

The same algebra $L$ provides a counter-example in characteristic $p$
to the statement of Theorem \ref{th:rad},
because it is solvable but not all of its elements
are $\overrightarrow{v}$-Engel (for example, $y$ as above is not).
\end{example}


Let us now consider another sequence $\overrightarrow{w}$ which
hopefully will be extendable to the case of positive
characteristic. Define $w_1(x,y)=[x,y]$ and, by induction,
$w_{n+1}(x,y)=[[w_n(x,y),x],[w_n(x,y),y]]$.

\begin{theorem} \cite{GKNP} \label{th:ch-w}
Let $L$ be a finite dimensional Lie algebra over a field $k$ of
characteristic different from $2,3,5$. Then $L$ is solvable if and
only if for some $n$ the identity $w_n(x,y)\equiv 0$ holds in $L$.
\end{theorem}

We hope that the same sequence $w_n$ allows one to describe the
solvable radical of $L$.


\begin{theorem} \label{th:rad-w}
Let $L$ be a finite dimensional Lie algebra over an algebraically
closed field $k$ of characteristic zero. Then its solvable radical
$R$ coincides with the set of all $\overrightarrow{w}$-Engel
elements of $L$.
\end{theorem}

\begin{proof}
Let us first prove that every element $y$ of $R$ is
$\overrightarrow{w}$-Engel. Since $R$ is an ideal of $L$, we have
$w_1(x,y)\in R$ for any $x\in L$. From the definition of $w_n$ it
follows that $w_2(x,y)\in [R,R]$, and, by induction, $w_n(x,y)$
belongs to the $(n-1)$th term of the derived series of $R$. Since $R$
is solvable, for some $n$ we have $w_n(x,y)=0.$

In the opposite direction, we have to prove that every
$\overrightarrow{w}$-Engel element lies in the solvable radical.
It is enough to prove that there are no nonzero
$\overrightarrow{w}$-Engel elements in a simple $k$-algebra $L$.

Let $L$ be a simple Lie algebra over $k$.
As in the proof of Theorem \ref{th:rad},
we choose a Cartan subalgebra $H\subset L$ and a Borel subalgebra $B\supset H$,
and we choose $e_\alpha\in L_\alpha$, $e_\alpha\neq 0$, for all $\alpha\in \Phi=R(L,H)$.
We write an element $y\in L$ as
$y=y^- +h+ y^+$, where
$y^-\in L^-,\ h\in H,\ y^+\in L^+$.
Write
$y^-=\sum_{\alpha\in\Phi^+}c_{-\alpha}(y)e_{-\alpha}$.
Let $\gamma$ denote the maximal root in $\Phi$.

Let $y\in L,\ y\neq 0$.
As in the proof of Theorem \ref{th:rad},
we may and shall assume that
$c:=c_{-\gamma}(y)\neq 0$. We shall find $x$ such that
$w_n(x,y)\neq 0$ for any $n$. Take $x=e_\gamma$.
 Denote $h_{\gam}=[e_{\gam},e_{-\gam}]$, it is a nonzero element of $H$.
We have $[h_\gamma,e_\gamma]=a e_\gamma$ with $a\neq 0$, then
$[h_\gamma, e_{-\gamma}]=-ae_{-\gamma}$. Since $\gamma$ is the
maximal root, we have
$$
w_1(x,y)=[x,y]=ch_{\gam}+y_1^+ \text{ where } y_1^+\in L^+.
$$
Hence we have
$$
[w_1(x,y),x]=[c h_{\gam}+y_1^+,e_{\gam}]=c[h_\gamma,e_\gamma]=ca
e_{\gam}
$$
(we have $[y_1^+,e_{\gam}]=0$ because $\gam$ is the maximal root).
Furthermore,
$$
c_{-\gamma}([w_1(x,y),y])=c_{-\gamma}([ch_\gamma,y])=c_{-\gamma}([ch_\gamma,ce_{-\gamma}])=-c^2
a,
$$
(here once again we use the assumption that $\gam$ is the maximal
root). As at the first step, we obtain
$$
w_2(x,y)=[[w_1(x,y),x],[w_1(x,y),y]]=-c^3a^2h_{\gam}+y_2^+, \text{
where } y_2^+\in L^+.
$$
By induction we conclude that for all $n$ we have
$$
w_n(x,y)=(-c^2a^2)^{n-1} ch_\gamma +y_n^+\text{ where }y_n^+\in L^+.
$$
Since $c\neq 0$ and $a\neq 0$, we see that  $w_n(x,y)\neq 0$ for all $n$.
Thus there are no
nonzero $\overrightarrow{w}$-Engel elements in $L$.
\end{proof}


\begin{remark}
To extend Theorem \ref{th:rad-w} to the case of any field $k$ of
characteristic zero, it is enough to show that any
$\overrightarrow{w}$-Engel element $y\in L$ remains
$\overrightarrow{w}$-Engel in $\bar L=L\otimes_k\bar k$. To do
that, it suffices to reverse the order of quantifiers in the
definition of a $\overrightarrow{w}$-Engel element, i.e. to show
that if $y\in L$ is a $\overrightarrow{w}$-Engel element, then
there exists $n$ (depending only on $y$) such that for all $x\in
L$ we have $w_n(x,y)=0$.
This would imply that $y$ remains an Engel element in $\bar L$.
Indeed, with such a choice of $n$ the relation
$w_n(x,y)=0$ is a polynomial identity in $L$ (with respect to the variable
$x$), and it remains true as an identity in $\bar L$ (cf. Lemma
\ref{lem:cl}).

In other words, it is enough to show that the sequence
$A_n(y)=\{x\in L: w_n(x,y)=0\}$, $n=1,2,\dots$, of subsets of $L$
stabilizes provided $y$ is a $\overrightarrow{w}$-Engel element.
In light of Theorem \ref{th:rad-w}, this is valid in the case when
$k$ is an algebraically closed field. We can also prove this fact
for any uncountable field $k$, see Proposition \ref{contfield}
below.
\end{remark}

Let $k$ be a field. We fix an algebraically closed field $\Omega$
containing $k$. Let $S$ be a set of polynomials in the polynomial
ring $k[X_1,\dots,X_m]$ in $m$ variables. By a $k$-closed set (or
algebraic $k$-set) in the affine space $\BA^m$ we mean the set $B$
of common zeros in $\Omega^m$ of such a set of polynomials $S$
(cf. \cite[Ch.~IX, \S\S 1,2]{L}). In particular, we say that
$B=\BA^m$ if $B=\Omega^m$. By $B(k)$ we denote the set of zeros of
$S$ in $k^m$, i.e. $B(k)=B\cap k^m$.

\renewcommand{\a}{{\alpha}}

\begin{prop}\label {contfield}
Let $k$ be an uncountable field, and $\BA ^m$ be an affine space
over $k.$ Consider a sequence of  $k$-closed sets $B_n\subset \BA^m$,
such that $B_1\subseteq B_2\subseteq \dots B_n\subseteq\dots
$ and $\bigcup _{1}^{\infty} B_i(k)=\BA ^m(k).$ Then there is $n_0$ such
that $B_{n_0}=\BA ^m.$
\end{prop}

\begin{proof}
We use induction on $m.$

{\em Step} 1. If $m=1$ then any $k$-closed set is either finite or equal to
$\BA ^1.$ Thus,

--- either  $B_i$ are all finite, and $\bigcup _{1}^{\infty} B_i(k)$
contains at most countable set of elements, which contradicts to
the assumption  $\bigcup _{1}^{\infty} B_i(k)=\BA ^m(k);$

--- or there is $n_0$ such that $B_{n_0}=\BA ^1.$

{\em Step} $m$. Assume that the claim is valid for all $m'<m.$ Let $x_1,
\dots, x_m$ be the coordinates in $\BA ^m$, and let $H_0$ denote the hyperplane
defined by the equation $x_1=0$.

Fix a point $ \a=(0, a_2,\dots,a_m)\in H_0(k)$ and a line
$$
L_{\a}=\{(t,  a_2,\dots,a_m)\ |\  t\in \Omega\}.
$$
Since $ \cup (B_i\cap L_{\a})(k)= (\cup B_i(k))\cap L_\a(k)  =L_{\a}(k)$,
there is $n(\a)$ such that
$B_{n(\a)}\cap L_{\a}=L_{\a}$ (see Step 1).

For a natural number $l$ define a set
$$D_l=\{\a\in H_0(k): n(\a)\leq l\}\subset H_0(k).$$
In other words, $\a\in D_l$ means that $ L_{\a}\subset B_l.$ Let
us show that $D_l=V_l(k)$ for some  $k$-closed subset $V_l$ of
$H_0$.

Let $B_l$ be defined in $\BA ^m$ by a set of polynomial equations
\begin{equation*}
\begin {cases}
F_{1,l}(x_1,\dots ,x_m)=\sum x_1^ip_{i,1,l}(x_2,\dots ,x_m)=0;\\
F_{2,l}(x_1,\dots ,x_m)=\sum x_1^ip_{i,2,l}(x_2,\dots ,x_m)=0;\\
\dots \\
F_{s_l,l}(x_1,\dots ,x_m)=\sum x_1^ip_{i,s_l,l}(x_2,\dots
,x_m)=0.\end{cases}
\end{equation*}
The condition $L_{\a}\subset B_l$ implies that
$$F_{j,l}(t,a_2,\dots ,a_m)=0$$
for any  $t\in\Omega$, i.e. $p_{i,j,l}(a_2,\dots ,a_m)=0$ for all
possible $i,j$. Thus
$$D_l=\{\a\in H_0(k):\  p_{i,j,l}(\a)=0\},$$
and clearly $D_l=V_l(k)$ for the $k$-closed subset  $V_l$  of $H_0$ defined
by these equations.

On the other hand, $\bigcup V_l(k)=\bigcup D_l= H_0(k).$ Thus, by
the induction hypothesis there exists $l_0$ such that $
V_{l_0}=H_0,$ and, consequently, $B_{l_0}=\BA^m.$
\end{proof}

\begin{corollary} \label{cor:rad-w}
Let $L$ be a finite dimensional Lie algebra over an uncountable
field $k$ of characteristic zero. Then its solvable radical
$R$ coincides with the set of all $\overrightarrow{w}$-Engel
elements of $L$.
\end{corollary}

\begin{remark} \label{rem:count} The Proposition \ref{contfield}
is not valid for countable fields.

Here is an example. Let $x_1,\dots, x_n, \dots$ be the countable
set of all the $k$-points of the affine space $\BA^m$ over a
countable field $k$ for $m>1$. Denote by $L_n$ the straight line
containing $x_n$ and the origin, and set $B_n=\bigcup_{1}^{n}
L_i.$ Then  $\bigcup _{1}^{\infty} B_n(k) =\BA^m(k).$

Nevertheless, our conjecture is that Theorem \ref{th:rad-w} is valid
for  any field of characteristic zero.
\end{remark}

\begin{remark} \label{rem:p}
Theorem \ref{th:rad-w} does not hold in positive characteristic
because simple Lie algebras may then contain nonzero
$\overrightarrow{w}$-Engel elements (see Example \ref{ex:W} below).
However, we hope that
the theorem remains true for ``classical'' Lie algebras (i.e. for
those coming from algebraic groups).
\end{remark}


\begin{example} \label{ex:W}
Let $L=W(1;1)$ be the Witt algebra defined over a field $k$ of
characteristic $p$. Recall (see, for example, \cite[4.2, p.~148]{SF})
that $L$ is of dimension $p$ with multiplication table defined on a basis
$\{e_{-1}, e_0, e_1,\dots , e_{p-2}\}$ as follows:
\begin{equation} \label{eq:W}
[e_i,e_j]=\left\{
\begin{array}{cl}
(j-i)e_{i+j} & {\textrm{ if }} -1\leq i+j \leq p-2, \\
0            & {\textrm{ otherwise}}.
\end{array}
\right.
\end{equation}
If $p>2$, the algebra $L$ is simple \cite[Thm.~2.4(1) on p.~149]{SF}.
However, if $p>3$, it contains nonzero $\overrightarrow{w}$-Engel elements.
Indeed, let $y=e_{p-2}$, and let $x=\alp_{-1}e_{-1}+\dots +\alp_{p-2}e_{p-2}$
be an arbitrary element of $L$. From formulas (\ref{eq:W}) it follows
that
$$
w_1(x,y)=\alp_{-1}e_{p-3}+\alp_0e_{p-2}.
$$
For $p>3$ this implies $[w_1(x,y),y]=0$ and hence $w_2(x,y)=0$. Thus
$y$ is a $\overrightarrow{w}$-Engel element, and the statement
of Theorem \ref{th:rad-w} does not hold for $L$.
\end{example}


In order to extend the characterization of the solvable radical
obtained in Theorem \ref{th:rad-w} to a broader class of algebras,
we suggest a little change in our original approach.

\begin{defn}
We say that an element $y\in L$ is {\em strictly Engel}, if it is
Engel and for any $x\in L$ the element $[x,y]$ is Engel
(see Definition \ref{def:engel} for the notion of an Engel element).
\end{defn}

For finite dimensional Lie algebras the following counterpart of
Baer's theorem holds.

\begin{theorem} \label{cl}
Let $L$ be a finite dimensional Lie algebra a field $k$, $\char (k)=0$.
The nilpotent radical $N$ of $L$ coincides with the set of all
strictly Engel elements of $L$.
\end{theorem}

\begin{proof}
By \cite[Thm. II.3.3]{J} every element $y$ of $N$
is an Engel element.
Since $N$ is an ideal, $[x,y]\in N$ for every $x\in L,\ y\in N$.
Thus every element of $N$ is strictly Engel.

Let us show that there are no nonzero strictly Engel elements
outside $N$. Let  $y$ be  a strictly Engel element of $L$.
Denote by $R$ the solvable radical of $L$.
First we prove that $y\in R$ .
Indeed, if $y\notin R$, then by Theorem \ref{th:rad}
there exists $x\in L$ such that $(\ad [x,y])^n x\neq 0$ for any $n$.
Hence $[x,y]$ is not Engel and therefore $y$ is not strictly Engel.

Thus $y\in R$. The set of Engel elements of $L$ which are
contained in $R$ coincides with  $N$,
see \cite[Ch.~I, \S5, n${}^\circ$5, Cor.~7 of Thm.~1]{Bou} (see also
\cite[Ch.~16, \S4)]{AS}. Thus $y\in N$.
\end{proof}

\begin{remark}\label{rem:3.15}
There is another approach allowing one to characterize the nilpotent
radical in a slightly different manner.

Let $L$ be a finite dimensional Lie algebra, and $L^+$ be the
corresponding vector space. Denote by $(\ad L)^*$ the subalgebra
of the associative algebra $\End (L^+)$ generated by the
linear operators of the form $\ad x$ where $x\in L$.
Let $y\in L$. Observe that the principal ideal $\langle y\rangle$
in the Lie algebra $L$ generated by $y$
consists of the elements $u(y)$, where $u\in (\ad L)^*$.

We say that an Engel element $y$ is {\it totally Engel} if $u(y)$
is an Engel element for every $u\in (\ad L)^*$.

The next theorem follows from the definitions.

\begin{theorem}
Let $L$ be a finite dimensional Lie algebra over a field $k$.
The nilpotent radical $N$ of $L$ coincides with
the set of all totally Engel elements of $L$.
\end{theorem}

\begin{proof}
An argument similar to that  in the proof of Theorem
\ref{cl}, using \cite[Theorem II.3.3]{J}, shows that
all elements of $N$ are totally Engel.

Conversely, if $y$ is totally Engel, then the ideal $\langle y\rangle$
consists of Engel elements.
By Engel's theorem, this ideal is nilpotent and therefore is contained in $N$.
In particular $y\in N$.
\end{proof}


\end{remark}

\begin{remark}
Observations from the previous remark are
also relevant to the case of the solvable radical and
$\overrightarrow{w}$-Engel elements. We could define strictly (resp. totally)
$\overrightarrow{w}$-Engel elements in the same way as it was done
for strictly (resp. totally) Engel elements. Moreover, to treat the case of
positive characteristic, one may be led to an even more
restrictive definition of totally $\overrightarrow{w}$-Engel
element requiring that it remains $\overrightarrow{w}$-Engel after
applying any (not necessary inner) derivation of $L$. In this
setting one can hope to get a characterization of the solvable
radical as the set of all totally $\overrightarrow{w}$-Engel
elements.
\end{remark}

\begin{remark} \label{rem:three}
Using yet another approach, one can redefine the notion of
strictly $\overrightarrow{w}$-Engel element using sequences of
three (or more) variables. For example, define
$r_n(x,y,z)=[z,[x,y],\dots ,[x,y]]$. Then  an element $y\in L$ is
strictly Engel if for every $x,z$ there exists $n=n(x,y,z)$ such that
$r_n(x,y,z)=0$ and $e_n(x,y)=0$.
\end{remark}

\section{Group case. Reduction theorem} \label{sec:gr}

Let $G$ be a finite group. We have
$[x,y]=xyx^{-1}y^{-1}=x\sigma_y(x^{-1})$ where $\sigma_y\in\Aut G$
takes $x$ to $yxy^{-1}$.

Let $\overrightarrow{u}=u_1(x,y),\ldots, u_n(x,y),\dots$ be a
sequence in $F_2=F_2(x,y)$. We want to define
$\overrightarrow{u}$-Engel-like automorphisms of $G$.

Let $G_1=G\leftthreetimes \Aut G$ be the semidirect product of $G$
and $\Aut G$. Recall that it can be viewed as the set of all pairs
of the form $(g,\sigma)$, where $g\in G$, $\sigma\in \Aut G$, with
multiplication $(g_1,\sigma_1)\centerdot
(g_2,\sigma_2)=(g_1\sigma_1(g_2),\sigma_1\sigma_2)$. The natural
embeddings $G\to G_1$ and $\Aut G\to G_1$, together with the above
formula for the group law in $G_1$, imply that the equality
$\sigma g\sigma^{-1}=\sigma(g)$ holds inside $G_1$, for any $g\in
G$ and $\sigma\in \Aut G$ (informally, an arbitrary automorphism
of $G$ becomes an inner automorphism inside a bigger group $G_1$).
In particular, this implies that $\sigma g\sigma^{-1}$ belongs to
$G$.

A $\overrightarrow{u}$-Engel-like automorphism $\sigma$ of $G$
should be defined in such a way that in the group
$G_1=G\leftthreetimes \Aut G$ it will be presented by a usual
$\overrightarrow{u}$-Engel-like element.

Let us take a correct sequence $\overrightarrow{u}(x,y)$ in
$F_2(x,y)$. Given a group $G$, consider a homomorphism $\mu\colon
F_2(x,y) \to G_1$. Denote $\mu(x)=g,$ $\mu(y)=\sigma$, where $g\in
G$, $\sigma\in \Aut G$, and define $u_n(g,\sigma )=\mu
(u_n(x,y))\in G_1$ (informally, we ``substitute'' $g$ instead of
$x$ and $\sigma$ instead of $y$).

\begin{lemma} \label{lem:sigma} Suppose that a correct sequence
$\overrightarrow{u}=\overrightarrow{u}(x,y)$
in $F_2(x,y)$ satisfies the  property:

(iii) For every $x,y\in F_2$ the word $u_n(x,y)$ belongs to the
minimal subgroup containing the $\sigma_y$-orbit of the element
$x$.

Then for every group $G\leftthreetimes \Aut G$, every $g\in G$,
and every $\sigma\in\Aut G$, the element $u_n(g,\sigma)$ belongs
to $G$.
\end{lemma}

\begin{proof}
According to condition (iii), every element $u_n(x,y)$ can be
represented as a product of elements $yxy^{-1}$, $y^kx^{-1}y^{-k}.$
Then $u_n(g,\sigma)$ is a product of elements $\sigma^k g^{\pm
1}\sigma^{-k}=\sigma^k(g^{\pm 1})$ and thus belongs to $G$.
\end{proof}

\begin{defn} \label{def:auto}
We call a sequence $\overrightarrow{u}$ satisfying the hypotheses
of Lemma $\ref{lem:sigma}$ autocorrect.
\end{defn}

We are now able to define $\overrightarrow{u}$-Engel
automorphisms:


\begin{defn}
Let $\overrightarrow{u}$ be an autocorrect sequence in $F_2$ and
let $G$ be a group. We say that an
automorphism $\sigma\in \Aut G$ is $\overrightarrow{u}$-Engel if
for any $g\in G$ there exists $n$ such that $u_n(g,\sigma)=1$.
\end{defn}

\begin{example} \label{ex:sigma}
For $\sigma\in \Aut(G)$ define $[x,\sigma]=x\sigma(x^{-1})$.

Consider the Engel sequence $\overrightarrow{e}$ defined by
$$
e_1(x,\sigma)=[x,\sigma],\ldots,
e_n(x,\sigma)=[e_{n-1}(x,\sigma),\sigma].
$$
It satisfies the hypotheses of Lemma \ref{lem:sigma}.
We say that an automorphism $\sigma$ of $G$ is Engel
if for any $x\in G$ there exists $n$ such that $e_n(x,\sigma)=1$.


Another example of a sequence satisfying the hypotheses of Lemma
\ref{lem:sigma} is
$\overrightarrow{w}=\overrightarrow{w}(x,\sigma),$ where
$$
w_1=[x,\sigma], \ldots, w_n=[[w_{n-1},x],[w_{n-1},\sigma]].
$$
\end{example}

Note that $\sigma=\sigma_y$ is a $\overrightarrow{u}$-Engel automorphism if
and only if $y$ is a $\overrightarrow{u}$-Engel element. Thus if
$G$ has no non-trivial $\overrightarrow{u}$-Engel automorphisms,
then $G$ has no non-trivial $\overrightarrow{u}$-Engel elements.
Let $A$ be the cyclic subgroup of $\Aut(G)$ generated by $\sigma$.
Denote $\widetilde G=G \leftthreetimes A$. If $\sigma$ is a
$\overrightarrow{u}$-Engel automorphism, then $\sigma$ is a
$\overrightarrow{u}$-Engel element in $\widetilde G$.

Let $G$ be a finite semisimple group, i.e. $R(G)=1$, and let
$\overrightarrow{u}$ be an autocorrect sequence. 
Consider two classes of groups:

\begin{itemize}
\item the class of semisimple groups with non-trivial
$\overrightarrow{u}$-Engel elements;

\item the class of semisimple groups with non-trivial
$\overrightarrow{u}$-Engel automorphisms.
\end{itemize}
Our aim is to show that for some $\overrightarrow{u}$ these
classes are empty.

Any finite semisimple group $G$ contains a unique maximal normal
centreless completely reducible (CR) subgroup (by definition, CR
means a direct product of finite non-abelian simple groups) which
is called the CR-radical of $G$ (see \cite[3.3.16]{Ro}). We call a
product of the isomorphic factors in the decomposition of the
$CR$-radical {\it an isotypic component} of $G$.

Denote the $CR$-radical of $G$ by $V=V(G)$. This is a
characteristic subgroup of $G$.

\begin{theorem} \label{th:CR}
Let $\overrightarrow{u}$ be a correct sequence, and let $G$ be a
semisimple group of smallest order having non-trivial
$\overrightarrow{u}$-Engel elements. Then $G$ has the following
structure:

(i) all non-trivial quotients of $G$ are solvable;

(ii) the CR-radical of $G$ contains only one isotypic component.
\end{theorem}

\begin{proof}
(i) Let $H$ be an arbitrary normal subgroup in $G$. Denote $\bar G =
G/H$. Let $R(\bar G)=H_1/H$ be the solvable radical of $\bar G$.
Then $\bar{\bar G}= \bar G/R(\bar G)$ is semisimple and of order
strictly smaller than the order of $G$. Hence, for a
$\overrightarrow{u}$-Engel element $g\in G$ we have $\bar{\bar
g}=1$, and therefore $\bar g\in H_1/H$. Thus $g\in H_1$. Since
$H_1$ is normal in  $G$, it is semisimple. The order of $H_1$ is
smaller than the order of $G$ and $g\in H_1$, therefore $g=1$.

Thus $G=H_1$, and $G/H=H_1/H$ is solvable.

(ii) Take $H=V(G)$, where $V(G)$ is the CR-radical of $G$. Suppose
$V(G)=G_1\times G_2$, where $G_1$ and $G_2$ are isotypic
components of $V(G)$, i.e. products  of isomorphic non-abelian
simple groups. Then they are normal subgroups of $G$. By (i), $G/G_1$
and $G/G_2$ are solvable, hence so is $G/(G_1\cap G_2)$.
Since the intersection is 1, $G$ is solvable.
The contradiction shows that $V(G)$ consists of one isotypic component.
\end{proof}

Observe that under hypotheses of Theorem \ref{th:CR}, the
CR-radical of $G$ coincides with the intersection of all its
proper normal subgroups.


Let us now consider a group $G$ which is minimal with respect to
the second property (to possess non-trivial
$\overrightarrow{u}$-Engel automorphisms). We need the following
auxiliary result.

\begin{prop} \label{prop:action}
Let $G$ be the group generated by $G_0$ and $g_0$, and let
$G_0\triangleleft G$.  Let $G_0$ be a finite semisimple group, and
$V_0$ its CR-radical. If $g_0$ acts trivially on $V_0$ then $g_0$
acts trivially on $G_0.$
\end{prop}

\begin{proof}
Denote $R(G)=R$, $\bar G=G/R$. We have $R\cap G_0=1$ and
$[R,G_0]=1$, $\bar G_0=G_0R/R \cong G_0$, $\bar V_0=V_0R/R \cong
V_0$.

Consider the action of $g_0$ on $G_0$ and on $\bar G_0.$ We have
$\overline {g_0gg_0^{-1}}=\bar g_0 \bar g \bar g_0^{-1}=g_0\bar g
g_0^{-1}$, $g\in G_0$. Thus, the actions of $g_0$ in $G_0$ and
$\bar G_0$ are isomorphic. $\bar V_0$ is the CR-radical of each of the groups
$\bar G_0$ and $\bar G$. Suppose that $\bar g_0$ belongs to the centralizer of
$\bar V_0$. Then $\bar  g_0 =\bar 1$, and $g_0\in R(G)$, $g_0$ and
$G_0$ commute. Hence $g_0$ acts trivially on $G_0$.
\end{proof}

Proposition \ref{prop:action} implies that if $\sigma$ is an
automorphism of a semisimple group $G$ acting trivially on the
CR-radical of $G$, then $\sigma$ is the identity automorphism of
$G$.

\begin{theorem} \label{th:min}
Let $\overrightarrow{u}$ be an autocorrect sequence, and let $G$ be a
semisimple group of smallest order having a non-trivial
$\overrightarrow{u}$-Engel automorphism. Then $G$ is a direct
product of isomorphic non-abelian simple groups.
\end{theorem}

\begin{proof}
Observe, first of all, that if $\sigma$ is a
$\overrightarrow{u}$-Engel automorphism of $G$ and $H$ is a
$\sigma$-invariant subgroup of $G$, then the corresponding
restricted automorphism $\sigma_1$ is a $\overrightarrow{u}$-Engel
automorphism of $H$. If now $H$ is a normal subgroup of $G$
invariant under a $\overrightarrow{u}$-Engel automorphism
$\sigma$, then $\bar \sigma$, the automorphism of the group $G/H$
induced by $\sigma$, is a $\overrightarrow{u}$-Engel automorphism
of $G/H$.

Let $G$ be a minimal semisimple group which has a non-trivial
$\overrightarrow{u}$-Engel automorphism $\sigma$. Let $V$ be the
CR-radical of $G$. If $V<G$ then the restriction of $\sigma$ to
$V$ induces a $\overrightarrow{u}$-Engel automorphism $\sigma_1$
of $V$ which must be trivial. Then according to Proposition
\ref{prop:action}, $\sigma$ is trivial on $G$. Contradiction. Thus
$V=G$.

Suppose $V=G$ consists of several isotypic components
$V_1,V_2,\ldots,V_k$. By the minimality hypothesis on $G$,
$\sigma$ acts trivially on each component, and $\sigma$ once again
must be trivial. The contradiction shows that $V$ consists of one isotypic
component.
\end{proof}

\begin{conj}\label{conj:prod}
There exists an autocorrect sequence $\overrightarrow{u}$ such that no
group $G$ of the form $G=\prod G_i$, where all $G_i$'s are
isomorphic simple non-abelian groups, has a non-trivial
$\overrightarrow{u}$-Engel automorphism.
\end{conj}

\begin{theorem} \label{th:red}
Conjecture $\ref{conj:radgr}$ is equivalent to Conjecture
$\ref{conj:prod}$.
\end{theorem}

\begin{proof}
First note that for any $G$ and any correct sequence
$\overrightarrow{u}$ all elements of $R=R(G)$ are
$\overrightarrow{u}$-Engel. Hence the assertion of Conjecture
\ref{conj:radgr} is equivalent to the following one: if $G$ is
semisimple, then it contains no non-trivial
$\overrightarrow{u}$-Engel elements.

1. Suppose that $\overrightarrow{u}$ is a sequence such that for
any group $G$ we have $R(G)=G(\overrightarrow{u})$. Let $G$ be any
finite semisimple group, and let $\sigma$ be a non-trivial
$\overrightarrow{u}$-Engel automorphism of $G$. Let us show that
$\sigma =1$. Let $\widetilde G=G \leftthreetimes \left<\sigma\right>$, and
consider $\bar G=\widetilde G/R(\widetilde G)$. This
group is semisimple, and $\bar \sigma$ is a
$\overrightarrow{u}$-Engel element of the semisimple group $\bar
G$. Then by hypothesis $\bar\sigma=1$, and $\sigma$ is trivial.
Thus the statement of Conjecture \ref{conj:prod} is fulfilled.

2. Suppose Conjecture \ref{conj:prod} is true, and take
$\overrightarrow{u}$ as in its statement. We want to prove that
there are no non-trivial $\overrightarrow{u}$-Engel elements in
any semisimple group. Assume the contrary, and let $G$ be a
semisimple group of smallest order containing a non-trivial
$\overrightarrow{u}$-Engel element $g$. Then $\sigma_g$ is a
non-trivial $\overrightarrow{u}$-Engel automorphism of $G$.
According to Theorem \ref{th:min}, $G=\prod G_i$, where all
$G_i$'s are isomorphic non-abelian simple groups. Conjecture
\ref{conj:prod} gives $\sigma_g=1$. Hence $g=1$.
\end{proof}

\begin{remark} \label{rem:conjseq}
We believe that Conjecture \ref{conj:prod} can be replaced with
the following, {\it a priori} weaker

\begin{conj}\label{conj:seq}
There exists an autocorrect sequence $\overrightarrow{u}$ such that no
finite simple non-abelian group $G$ has a non-trivial
$\overrightarrow{u}$-Engel automorphism.
\end{conj}

Indeed, Theorem 3.3.20 of \cite{Ro} says that if $G=\prod_{i}
G_i$, where all $G_i$'s are isomorphic simple non-abelian groups,
then $\Aut G$ is isomorphic to the wreath product  $\Aut G_i\ wr\
S_n$, where $S_n$ is the symmetric group of degree $n$. This gives
some hope to deduce Conjecture \ref{conj:prod} from Conjecture
\ref{conj:seq}.
\end{remark}

\begin{remark}
Conjectures \ref{conj:seq} and \ref{conj:prod} require to specify
a sequence $\overrightarrow{u}$. One of possible candidates is
the sequence $\overrightarrow{w}$ from Example \ref{ex:sigma}.
\end{remark}



\end{document}